\newif\iffinal
\finalfalse % Not final version
\finaltrue % Final version

\documentclass[reqno,twoside,11pt]{amsart}
\iffinal\else\usepackage[notref,notcite]{showkeys}\fi
\usepackage{amsmath}
\usepackage{amsfonts}
\usepackage{amssymb}
\usepackage{verbatim}
\usepackage{epsfig}

\IfFileExists{epsf.def}{\input epsf.def}{\usepackage{epsf}}

\IfFileExists{myowntimes.sty}{\usepackage{myowntimes}\usepackage{mathrsfs}}
    {\usepackage{times}\usepackage{mathrsfs}}
\DeclareFontFamily{OT1}{eusb}{} \DeclareFontShape{OT1}{eusb}{m}{n} {<5> <6> <7> <8> <9> <10> <11> <12> <14.4> eusb10}{}
\DeclareMathAlphabet{\eusb}{OT1}{eusb}{m}{n}

\DeclareFontFamily{OT1}{eusm}{} \DeclareFontShape{OT1}{eusm}{m}{n} {<5> <6> <7> <8> <9> <10> <11> <12> <14.4> eusm10}{}
\DeclareMathAlphabet{\eusm}{OT1}{eusm}{m}{n}

\DeclareFontFamily{OT1}{eufm}{} \DeclareFontShape{OT1}{eufm}{m}{n} {<5> <6> <7> <8> <9> <10> <11> <12> <14.4> eufm10}{}
\DeclareMathAlphabet{\mathfrak}{OT1}{eufm}{m}{n}

\DeclareFontFamily{OT1}{fraktura}{}
\DeclareFontShape{OT1}{fraktura}{m}{n} {<5> <6> <7> <8> <9> <10> <11> <12> <13> <14.4> [1.1] eufm10}{}
\DeclareMathAlphabet{\fraktura}{OT1}{fraktura}{m}{n}

\DeclareFontFamily{OT1}{cmfi}{} \DeclareFontShape{OT1}{cmfi}{m}{n} {<5> <6> <7> <8> <9> <10> <11> <12> <13> <14.4> [0.9] cmfi10}{}
\DeclareMathAlphabet{\cmfi}{OT1}{cmfi}{b}{n}

\DeclareFontFamily{OT1}{cmss}{} \DeclareFontShape{OT1}{cmss}{m}{n} {<5> <6> <7> <8> <9> <10> <11> <12> <13> <14.4> cmss10}{}
\DeclareMathAlphabet{\cmss}{OT1}{cmss}{m}{n}
\newtheoremstyle{thm}{1.5ex}{1.5ex}{\itshape\rmfamily}{} {\bfseries\rmfamily}{}{2ex}{}

\newtheoremstyle{def}{1.5ex}{1.5ex}{\slshape\rmfamily}{} {\bfseries\rmfamily}{}{2ex}{}

\newtheoremstyle{rem}{1.3ex}{1.3ex}{\rmfamily}{} {\itshape}
{} {1.5ex}{}

\newenvironment{proofsect}[1] {\vskip0.1cm\noindent{\rmfamily\itshape#1.}}{\qed\vspace{0.15cm}}%{\newline\vspace{0.15cm}}

\theoremstyle{thm}
\newtheorem{theorem}{Theorem}[section]
\newtheorem{lemma}[theorem]{Lemma}
\newtheorem{proposition}[theorem]{Proposition}

\newtheorem*{Main Theorem}{Main Theorem.}

\newtheorem*{special theorem}{Lindeberg-Feller Theorem for Martingales}

\theoremstyle{def}

\theoremstyle{rem}

\numberwithin{equation}{section}

\renewcommand{\section}{\secdef\sct\sect}
\newcommand{\sct}[2][default]{%
\refstepcounter{section}
\addcontentsline{toc}{section}{{\tocsection {}{\thesection}{\!\!\!\!#1\dotfill}}{}}
\vspace{0.7cm}
\centerline{\scshape\thesection.\ #1} \nopagebreak \vspace{0.2cm}}
\newcommand{\sect}[1]{%
\vspace{0.4cm} \centerline{\large\scshape\rmfamily #1}
\vspace{0.2cm}}

\renewcommand{\subsection}{\secdef\subsct\sbsect}
\newcommand{\subsct}[2][default]{\refstepcounter{subsection}
\addcontentsline{toc}{subsection}
{{\tocsection{\!\!}{\hspace{1.2em}\thesubsection}{\!\!\!\!#1\dotfill}}{}}
\nopagebreak\vspace{0.45\baselineskip} {\flushleft\bf
\thesubsection~\bf #1.~}
\\*[3mm]\noindent
\nopagebreak}
\newcommand{\sbsect}[1]{\vspace{0.1cm}\noindent
\textbf{#1.~}\vspace{0.1cm}}

\renewcommand{\subsubsection}{%
\secdef \subsubsect\sbsbsect}
\newcommand{\subsubsect}[2][default]{%
\refstepcounter{subsubsection}
\addcontentsline{toc}{subsubsection}{{\tocsection{\!\!}
{\hspace{3.05em}\thesubsubsection}{\!\!\!\!#1\dotfill}}{}}
\nopagebreak
\vspace{0.15\baselineskip} \nopagebreak {\flushleft\rmfamily
\itshape\thesubsubsection
\ \rmfamily #1\/.}\ }
\newcommand{\sbsbsect}[1]{\vspace{0.1cm}\noindent
\rmfamily \itshape
\arabic{section}.\arabic{subsection}.\arabic{subsubsection} \
\sffamily #1\/.\ }

\iffinal

\else

\fi

\renewcommand{\caption}[1]{%
\vglue0.5cm
\refstepcounter{figure}
\begin{minipage}{0.9\textwidth}\small {\sc Figure~\thefigure. }#1\end{minipage}}

\newcommand{\diam}{\operatorname{diam}}

\newcommand{\textd}{\text{\rm d}\mkern0.5mu}

\newcommand{\texte}{\text{\rm e}}

\newcommand{\1}{\operatorname{\sf 1}}

\newcommand{\FF}{\mathcal F}
\newcommand{\GG}{\mathcal G}
\newcommand{\HH}{\mathcal H}

\newcommand{\LL}{\mathcal L}

\newcommand{\B}{\mathbb B}

\newcommand{\E}{\mathbb E}

\newcommand{\G}{\mathbb G}

\newcommand{\BbbP}{\mathbb P}
\newcommand{\Q}{\mathbb Q}

\newcommand{\T}{\mathbb T}

\newcommand{\Z}{\mathbb Z}

\newcommand{\scrC}{\mathscr{C}}

\newcommand{\twoeqref}[2]{(\ref{#1}--\ref{#2})}
\newcommand{\cc}{{\text{\rm c}}}

\def\myffrac#1#2 in #3{\raise 2.6pt\hbox{$#3 #1$}\mkern-1.5mu\raise 0.8pt\hbox{$#3/$}\mkern-1.1mu\lower 1.5pt\hbox{$#3 #2$}}
\newcommand{\ffrac}[2]{\mathchoice%
    {\myffrac{#1}{#2} in \scriptstyle}
    {\myffrac{#1}{#2} in \scriptstyle}
    {\myffrac{#1}{#2} in \scriptscriptstyle}
    {\myffrac{#1}{#2} in \scriptscriptstyle}
}

\newcommand{\hate}{\hat{\text{\rm e}}}

\newcommand{\pomega}{\partial^{\,\omega}\!}
\newcommand{\pomegaalpha}{\partial^{\,\omega,\alpha}\!}
\newcommand{\pstar}{\partial^{\,*}\!}

\title[Random walk among random conductances]
{\fontsize{13}{16}\selectfont Anomalous heat-kernel decay for random walk among bounded random conductances}
\author[N.~Berger, M.~Biskup, C.E.~Hoffman and G.~Kozma]{N.~Berger\,$^1$,\, M.~Biskup\,$^1$,\, C.E.~Hoffman\,$^2$ \and\, G.~Kozma\,$^3$}

\begin{document}
\thanks{\hglue-4.5mm\fontsize{9.6}{9.6}\selectfont\copyright\,2007 N.~Berger, M.~Biskup, C.~Hoffman and G.~Kozma. Reproduction, by any means, of the entire
article for non-commercial purposes is permitted without charge.\vspace{2mm}}
\maketitle

\vspace{-4mm}
\centerline{\textit{$^1$Department of Mathematics, UCLA, Los Angeles, U.S.A.}}
\centerline{\textit{$^2$Department of Mathematics, University of Washington, Seattle, U.S.A.}}
\centerline{\textit{$^3$Weizmann Institute of Science, Rehovot, Israel}}

\vspace{-2mm}
\begin{abstract}
We consider the nearest-neighbor simple random walk on~$\Z^d$, $d\ge2$, driven by a field of bounded random conductances $\omega_{xy}\in[0,1]$. The conductance law is i.i.d.\ subject to the condition that the probability of $\omega_{xy}>0$ exceeds the threshold for bond percolation on~$\Z^d$. For environments in which the origin is connected to infinity by bonds with positive conductances, we study the decay of the $2n$-step return probability $\cmss P_\omega^{2n}(0,0)$. We prove that $\cmss P_\omega^{2n}(0,0)$ is bounded by a random constant times $n^{-d/2}$ in $d=2,3$, while it is $o(n^{-2})$ in~$d\ge5$ and $O(n^{-2}\log n)$ in $d=4$. By producing examples with anomalous heat-kernel decay approaching~$1/n^2$ we prove that the $o(n^{-2})$ bound in $d\ge5$ is the best possible. We also construct natural $n$-dependent environments that exhibit the extra $\log n$ factor in $d=4$. 
\end{abstract}

\section{Introduction}
\label{sec1}\noindent
Random walk in reversible random environments is one of the best studied subfields of random motion in random media. In continuous time, such walks are usually defined by their generators~$\LL_\omega$ which are of the form
\begin{equation}
%\label{}
(\LL_\omega f)(x)=\sum_{y\in\Z^d}\omega_{xy}\bigl[f(y)-f(x)\bigr],
\end{equation}
where $(\omega_{xy})$ is a family of random (non-negative) conductances subject to the symmetry condition~$\omega_{xy}=\omega_{yx}$. The sum $\pi_\omega(x)=\sum_y\omega_{xy}$ defines an invariant, reversible measure for the corresponding continuous-time Markov chain. The discrete-time walk shares the same reversible measure and is driven by the transition matrix
\begin{equation}
%\label{}
\cmss P_\omega(x,y)=\frac{\omega_{xy}}{\pi_\omega(x)}.
\end{equation}
In most situations~$\omega_{xy}$ are non-zero only for nearest neighbors on~$\Z^d$ and are sampled from a shift-invariant, ergodic or even i.i.d.\ measure~$\BbbP$ (with expectation henceforth denoted by~$\E$).

Two general classes of results are available for such random walks under the additional assumptions of uniform ellipticity,
\begin{equation}
\label{ellipticity}
\exists\alpha>0\,\colon\quad \BbbP(\alpha<\omega_b<\ffrac1\alpha)=1,
\end{equation}
and the boundedness of the jump distribution,
\begin{equation}
%\label{}
\exists R<\infty\,\colon\quad |x|\ge R\quad\Rightarrow\quad\cmss P_\omega(0,x)=0,\quad \BbbP\text{\rm-a.s.}
\end{equation}
First, as proved by Delmotte~\cite{Delmotte}, one has the standard, local-CLT like decay of the heat kernel ($c_1,c_2$ are absolute constants):
\begin{equation}
\label{heat-kernel}
\cmss P_\omega^n(x,y)\le \frac{c_1}{n^{d/2}}\exp\Bigl\{-c_2\frac{|x-y|^2}n\Bigr\}.
\end{equation}
Second, an annealed invariance principle holds in the sense that the law of the paths under the measure integrated over the environment scales to a non-degenerate Brownian motion (Kipnis and Varadhan~\cite{Kipnis-Varadhan}). A quenched invariance principle can also be proved by invoking techniques of homogenization theory (Sidoravicius and Sznitman~\cite{Sidoravicius-Sznitman}).

Once the assumption of uniform ellipticity is relaxed, matters get more complicated. The most-intensely studied example is the simple random walk on the infinite cluster of supercritical bond percolation on~$\Z^d$, $d\ge2$. This corresponds to~$\omega_{xy}\in\{0,1\}$ i.i.d.\ with~$\BbbP(\omega_b=1)>p_\cc(d)$ where~$p_\cc(d)$ is the percolation threshold. Here an annealed invariance principle has been obtained by De Masi, Ferrari, Goldstein and Wick~\cite{demas1,demas2} in the late 1980s. More recently, Mathieu and Remy~\cite{Mathieu-Remy} proved the on-diagonal (i.e., $x=y$) version of the heat-kernel upper bound \eqref{heat-kernel}---a slightly weaker version of which was also obtained by Heicklen and Hoffman~\cite{Heicklen-Hoffman}---and, soon afterwards, Barlow~\cite{Barlow} proved the full upper and lower bounds on $\cmss P_\omega^n(x,y)$ of the form \eqref{heat-kernel}. (Both of these results hold for $n$ exceeding some random time defined relative to the environment in the vicinity of~$x$ and~$y$.) Heat-kernel upper bounds were then used in the proofs of quenched invariance principles by Sidoravicius and Sznitman~\cite{Sidoravicius-Sznitman} for $d\ge4$, and for all $d\ge2$ by Berger and Biskup~\cite{BB} and Mathieu and Piatnitski~\cite{Mathieu-Piatnitski}.

\smallskip
Notwithstanding our precise definition \eqref{ellipticity}, the case of supercritical percolation may still be regarded as uniformly elliptic because the conductances on the percolation cluster are still uniformly bounded away from zero and infinity. It is thus not clear what phenomena we might encounter if we relax the uniform ellipticity assumption in an essential way. A number of quantities are expected (or can be proved) to vary continuously with the conductance distribution, e.g., the diffusive constant of the limiting Brownian motion. However, this may not apply to asymptotic statements like the heat-kernel bound \eqref{heat-kernel}.

In a recent paper, Fontes and Mathieu~\cite{Fontes-Mathieu} studied continuous-time random walk on~$\Z^d$ with conductances given by
\begin{equation}
\label{omega-wedge}
\omega_{xy}=\omega(x)\wedge\omega(y)
\end{equation}
for some i.i.d.\ random variables~$\omega(x)>0$. For these cases
it was found that the annealed heat  kernel,
$\E[P_{\omega,0}(X_t=0)]$, where~$P_{\omega,0}$ is the law of the
walk started at the origin and~$\E$ is the expectation with
respect to the environment, exhibits an \emph{anomalous decay} for
environments with too heavy lower tails at zero. Explicitly, from \cite[Theorem~4.3]{Fontes-Mathieu} we have
\begin{equation}
\label{annealed-HC}
\E\bigl[P_{\omega,0}(X_t=0)\bigr]=t^{-(\gamma\wedge\frac
d2)+o(1)},\qquad t\to\infty,
\end{equation}
where~$\gamma>0$ characterizes the lower tail of the~$\omega$-variables,
\begin{equation}
\label{MF-tail}
\BbbP\bigl(\omega(x)\le s\bigr)\sim s^\gamma,\qquad s\downarrow0.
\end{equation}
As for the quenched problem, for~$\gamma<\ffrac
d2$,~\cite[Theorem~5.1]{Fontes-Mathieu} provides a lower bound on the diagonal
heat-kernel decay exponent (a.k.a.\ spectral dimension):
\begin{equation}
%\label{}
\BbbP\bigl[P_{\omega,0}(X_t=0)\le
t^{-\alpha}\bigr]\,\underset{t\to\infty}\longrightarrow\,1
\end{equation}
for every~$\alpha<\alpha_0$ where
\begin{equation}
%\label{}
\alpha_0=\frac d2\frac{1+\gamma}{1+\ffrac d2}.
\end{equation}
But, since~$\alpha_0<\ffrac d2$, this does not rule out the usual diffusive scaling. Nevertheless, as~$\alpha_0>\gamma$ for~$\gamma<\ffrac d2$, the annealed and quenched heat-kernel decay at different rates.

\smallskip
The reason why the annealed heat  kernel may decay slower than usual can be seen rather directly from the following argument: The quenched probability that the walk does not even leave the origin up to time~$t$ is $\texte^{-t\pi_\omega(0)}$. By~$\pi_\omega(0)\le2d\omega(0)$, we have
\begin{equation}
\E\bigl[P_{\omega,0}(X_t=0)\bigr]\ge \E\texte^{-2d\omega(0) t}.
\end{equation}
For $\omega(0)$ with the tail \eqref{MF-tail}, this yields a lower bound of~$t^{-\gamma}$. (A deeper analysis shows that this is actually a dominating strategy~\cite{Fontes-Mathieu}.) A similar phenomenon can clearly be induced for $\omega_{xy}$ that are i.i.d.\ with a sufficiently heavy tail at zero, even though then the correspondence of the exponents in \twoeqref{annealed-HC}{MF-tail} will take a slightly different form.

The fact that the dominating strategy is so simple makes one wonder how much of this phenomenon is simply an artifact of taking the annealed average. Of not much help in this matter is the main result (Theorem~3.3) of Fontes and Mathieu~\cite{Fontes-Mathieu} which shows that the mixing time for the random walk on the largest connected component of a torus will exhibit anomalous (quenched) decay once $\gamma<\ffrac d2$. Indeed, the mixing time is by definition dominated by the worst-case local configurations that one can find \emph{anywhere} on the torus and thus the reasoning we used to explain the anomalous decay of the annealed heat  kernel applies here as well.

\smallskip
The main goal of this paper is to provide universal upper bounds on the \emph{quenched} heat kernel and support them by examples exhibiting the corresponding lower bounds. Somewhat surprisingly, and unlike for the annealed heat kernel, the existence of anomalous quenched heat-kernel decay turns out to be dimension dependent.

\section{Main results}
\noindent
We will work with a collection of bounded, nearest-neighbor conductances~$(\omega_b)\in\Omega=[0,1]^{\B}$ where~$b$ ranges over the set~$\B$ of unordered pairs of nearest neighbors in~$\Z^d$. The law $\BbbP$ of the~$\omega$'s will be i.i.d.\ subject to the condition that the bonds with positive conductances percolate. Given~$\omega$, we use~$\scrC_\infty=\scrC_\infty(\omega)$ to denote the set of sites that have a path to infinity along bonds with positive conductances.
It is well known that~$\scrC_\infty$ is connected with probability one.

\smallskip
The main result of this paper is as follows:

\begin{theorem}
\label{thm1}
Let~$d\ge2$ and consider a collection $\omega=(\omega_b)$ of i.i.d.\ conductances in $[0,1]$ with $\BbbP(\omega_b>0)>p_\cc(d)$ where~$p_\cc(d)$ is the threshold for bond percolation on~$\Z^d$. For almost every~$\omega\in\{0\in\scrC_\infty\}$, there is~$C=C(\omega)<\infty$ such that
\begin{equation}
\cmss P_\omega^n(0,0)\le C(\omega)\,\begin{cases}
n^{-d/2},\qquad&d=2,3,
\\
n^{-2}\log n,\qquad&d=4,
\\
n^{-2},\qquad&d\ge5,
\end{cases}
\end{equation}
for all~$n\ge1$. In fact, for $d\geq 5$, almost surely
\begin{equation}
\label{eq:eqnum1}
\lim_{n\to\infty}\,n^2\,\cmss P_\omega^n(0,0) = 0.
\end{equation}
\end{theorem}

Note that these estimates imply that the random walk is almost surely transient in all dimensions~$d\ge3$. This is of course a consequence of the fact---to be exploited in more depth later---that under~$p>p_\cc(d)$ one has an infinite cluster of bonds with conductances bounded strictly from below. Then a.s.\ transience in~$d\ge3$ follows by monotonicity in conductances and the result of Grimmett, Kesten and Zhang~\cite{Grimmett-Kesten-Zhang}. (Recurrence in~$d=1,2$ is inferred directly from the monotonicity of this notion in the conductances.)

To show that our general upper bound in~$d\ge5$ represents a real phenomenon, we state the existence of appropriate examples:

\begin{theorem}
\label{thm2}
(1)
Let~$d\ge5$ and $\kappa>\ffrac1d$. There exists an i.i.d.\
law~$\BbbP$ on bounded, nearest-neighbor conductances
with~$\BbbP(\omega_b>0)>p_\cc(d)$ and a random
variable~$C=C(\omega)$ such that for almost
every~$\omega\in\{0\in\scrC_\infty\}$,
\begin{equation}
\label{lower-bd}
\cmss P_\omega^{2n}(0,0)\ge\ C(\omega)\frac{\texte^{-(\log n)^\kappa}}{n^2},
\qquad n\ge1.
\end{equation}

\noindent
(2) Let $d\ge5$. For every increasing sequence
$\{\lambda_n\}_{n=1}^\infty$, $\lambda_n\to\infty$, there exists an i.i.d.\ law $\BbbP$ on
bounded, nearest-neighbor conductances
with~$\BbbP(\omega_b>0)>p_\cc(d)$ and an a.s.\ positive random variable~$C=C(\omega)$ such that for almost
every~$\omega\in\{0\in\scrC_\infty\}$,
\begin{equation}
\label{2.4}
\cmss P_\omega^n(0,0)\ge \frac{C(\omega)}{\lambda_nn^2}
\end{equation}
along a subsequence that does not depend on~$\omega$.
\end{theorem}

The upper bounds in Theorem~\ref{thm1} can be extended to more general
shift-invariant, ergodic environments under suitable assumptions
on their percolation properties. In particular, it follows that
for the Fontes-Mathieu example \twoeqref{annealed-HC}{MF-tail} no
anomaly occurs for the quenched heat  kernel in dimensions~$d=2,3$.
On the other hand, Theorem~\ref{thm2} can be specialized to the
case \eqref{omega-wedge} with i.i.d.~$\omega(x)$'s and, in~$d\ge5$, we can
produce anomalous decay as soon as the tails of $\omega$ at zero are
sufficiently heavy. (The constructions in the proof of Theorem~\ref{thm2} actually work for all $d\ge2$ but the result is of course interesting only for $d\ge5$.)

The distributions that we use in part~(1) of Theorem~\ref{thm2} have a tail near zero of the general form
\begin{equation}
\BbbP(\omega_{xy}<s) \approx |\log(s)|^{-\theta}
\end{equation}
with~$\theta>0$. Presumably, one can come up with examples of distributions that exhibit ``anomalous" behavior and have the power law tail,
\begin{equation}
\BbbP(\omega_{xy}<s) \approx s^{\gamma},
\end{equation}
for some $\gamma>0$. However, the construction seems to require subtle control of heat-kernel \emph{lower} bounds which go beyond the estimates that can be easily pulled out from the literature.

\smallskip
As we will see in the proofs, the underlying idea of all examples in Theorem~\ref{thm2} is the same: The walk finds a \emph{trap} which, by our specific choice, is a ``strong'' edge that can be reached only by crossing an edge of strength of order~$\ffrac1n$. Such traps allow the walk to get stuck for time of order~$n$ and thus improve its chances to make it back to the origin at the required time. To enter and exit the trap, the walk has to make two steps over the $O(\ffrac1n)$-edge; these are responsible for the overall $n^{-2}$-decay. Of course, in $d=2,3$ this cannot compete with the ``usual'' decay of the heat kernel and so we have to go to~$d\ge4$ to make this strategy dominant.

The upper bound in \eqref{eq:eqnum1} and the lower bound in \eqref{2.4} show that the $1/n^2$ decay in~$d\ge5$ is never achieved, but can be approached arbitrary closely. We believe the same holds also for $d=4$ for the decay rate $n^{-2}\log n$. We demonstrate the reason for our optimism by proving a lower bound for environments where the aforementioned traps occur with a positive density:

\begin{theorem}
\label{thm3}
Let~$d\ge4$ and let~$p>p_\cc(d)$. Sample a percolation configuration $\tilde\omega$ with parameter~$p$. For each~$n\ge1$ consider the i.i.d.\ environment $\omega^{(n)}$ defined from $\tilde\omega$ by putting $\omega^{(n)}_b=1$ on occupied bonds and~$\omega^{(n)}_b=\ffrac1n$ on vacant bonds. For a.e.~$\tilde\omega$ in which~$0$ has an occupied path to infinity, there is $C(\tilde\omega)>0$ such that for all~$n\ge1$,
\begin{equation}
\label{eq2.7}
\cmss P^{2n}_{\omega^{(n)}}(0,0)\ge C(\tilde\omega)
\begin{cases}
n^{-2},\quad&d\ge5,
\\
n^{-2}\log n,\quad&d=4.
\end{cases}
\end{equation}
\end{theorem}

We conclude with a remark concerning the path properties of the
above random walk. As mentioned previously, heat-kernel estimates
of the form~\eqref{heat-kernel} have been crucial for the proof of
the quenched invariance principle for simple random walk on
supercritical percolation clusters in~$d\ge3$. (The~$d=2$ argument
of Berger and Biskup~\cite{BB} actually avoids these bounds by appealing to the nearest-neighbor structure of the walk and to an underlying maximum principle.) The
absence of ``usual'' decay might suggest difficulty in following
the same strategy. Notwithstanding, using truncation to a ``strong component,'' a
version of which is invoked also in the present paper, this
problem can be circumvented and the corresponding quenched
invariance principle proved~(Mathieu~\cite{Mathieu-CLT}, Biskup and Prescott~\cite{BP}). 

Thus there are i.i.d.\ environments for which one has a functional CLT \emph{without} a local~CLT. This should not be too surprising as a CLT describes the
typical behavior whereas the heat-kernel decay, and a local-CLT,
describe rare events. Naturally, a CLT is much more robust
than its local counterpart.

\smallskip
Theorem~\ref{thm1} is proved in Sect.~\ref{sec3} while Theorems~\ref{thm2}-\ref{thm3} are proved in Sect.~\ref{sec4}. The Appendix (Sect.~\ref{sec5}) contains a self-contained proof of the isoperimetric inequality on the supercritical percolation cluster that we need in the proof of Theorem~\ref{thm1}.

\section{Heat-kernel upper bounds}
\label{sec3}\noindent
Here we will prove the heat-kernel bounds from Theorem~\ref{thm1}. The general strategy of our proof is as follows: For every~$\alpha>0$, we use~$\scrC_{\infty,\alpha}=\scrC_{\infty,\alpha}(\omega)$ to denote the set of all sites in~$\Z^d$ that have a path to infinity along edges with conductances at least~$\alpha$. Clearly, $\scrC_{\infty,\alpha}$ is a subgraph of~$\scrC_\infty$; we will sometimes refer to~$\scrC_{\infty,\alpha}$ as the \emph{strong component}.
We first prove the ``standard'' heat-kernel decay for the Markov chain obtained by recording the position of the random walk when it is on the strong component~$\scrC_{\infty,\alpha}$ for an appropriately chosen~$\alpha$. Then we control the difference between the time scales for the two walks using rather straightforward estimates.

\subsection{Coarse-grained walk}
The i.i.d.\ nature of the measure~$\BbbP$ ensures there is an a.s.-unique infinite connected component~$\scrC_\infty$ of bonds with positive conductances. Given $z\in\scrC_\infty$, we define the random walk~$X=(X_n)$ as a Markov chain on~$\scrC_\infty$ with transition probabilities
\begin{equation}
%\label{}
P_{\omega,z}(X_{n+1}=y|X_n=x)=\cmss P_\omega(x,y)=\frac{\omega_{xy}}{\pi_\omega(x)}
\end{equation}
and initial condition
\begin{equation}
%\label{}
P_{\omega,z}(X_0=z)=1.
\end{equation}
We use~$E_{\omega,z}$ to denote expectation with respect to~$P_{\omega,z}$. (Note the typographical distinction between the path distribution~$P_{\omega,z}$, the heat kernel~$\cmss P_\omega$, and the law of the environment~$\BbbP$.)

Next we will disregard bonds whose conductance is less than some small positive number~$\alpha$ which is chosen so that the remaining bonds still form an infinite component---to be denoted by~$\scrC_{\infty,\alpha}$. We quote Proposition~2.2 from Biskup and Prescott~\cite{BP}:

\begin{lemma}
\label{prop-percolace}
Let $d\ge2$ and~$p=\BbbP(\omega_b>0)>p_\cc(d)$. Then there exists $c(p,d)>0$ such that if~$\alpha$ satisfies
\begin{equation}
\label{infinite}
\BbbP(\omega_b\ge\alpha)>p_\cc(d)
\end{equation}
and
\begin{equation}
\label{finite}
\BbbP(0<\omega_b<\alpha)<c(p,d)
\end{equation}
then $\scrC_{\infty,\alpha}$ is nonempty and $\scrC_\infty\setminus\scrC_{\infty,\alpha}$ has only finite components a.s. In fact, if $\FF_x$ is the set of sites (possibly  empty) in the finite component of $\scrC_\infty\setminus\scrC_{\infty,\alpha}$ containing~$x$, then
\begin{equation}
\label{exponential}
\BbbP\bigl(x\in\scrC_\infty\,\,\&\,\,\diam \FF_x\ge n\bigr)\le C\texte^{-\eta n},\qquad n\ge1,
\end{equation}
for some~$C<\infty$ and $\eta>0$. Here ``$\diam$'' is the diameter in the~$\ell_\infty$-distance on~$\Z^d$.
\end{lemma}

Given~$z\in\scrC_{\infty,\alpha}$ we consider the \emph{coarse grained} random walk $\hat X=(\hat X_\ell)$---started at~$z$---which records the successive visits of~$X=(X_n)$ to $\scrC_{\infty,\alpha}$. Explicitly, let $T_1,T_2,\dots$ denote the times~$X$ takes between the successive steps of~$\hat X$, i.e., $T_{\ell+1}=\inf\{n>0\colon X_{T_\ell+n}\in\scrC_{\infty,\alpha}\}$ with $T_0=0$. Note that, as all components of $\scrC_\infty\setminus\scrC_{\infty,\alpha}$ are finite, $T_\ell<\infty$ a.s.\ for all~$\ell$. Then
\begin{equation}
%\label{}
\hat X_\ell=X_{T_1+\dots+T_\ell},\qquad \ell\ge1.
\end{equation}
Let $\hat{\cmss P}_\omega(x,y)$ denote the transition probability of the random walk $\hat X$,
\begin{equation}
%\label{}
\hat{\cmss P}_\omega(x,y)=P_{\omega,x}(X_{T_1}=y),\qquad x,y\in\scrC_{\infty,\alpha}.
\end{equation}
As is easy to check, the restriction of the measure $\pi_\omega$ to $\scrC_{\infty,\alpha}$ is invariant and reversible for the Markov chain on~$\scrC_{\infty,\alpha}$ induced by $\hat{\cmss P}_\omega$.

Consider the quantities
\begin{equation}
%\label{}
\hat\omega_{xy}=\pi_\omega(x)\hat{\cmss P}_\omega(x,y),\qquad x,y\in\scrC_{\infty,\alpha}.
\end{equation}
We may think of~$\hat X$ as the walk on~$\scrC_{\infty,\alpha}$
with the weak components ``re-wired'' by putting a bond with
conductance~$\hat\omega_{xy}$ between any pair of sites $(x,y)$ on
their (strong) boundary. By Lemma~\ref{prop-percolace}, all weak components are
finite and everything is well defined.

\smallskip
Our first item of business is to show that $\hat X$ obeys the standard heat-kernel bound:

\begin{lemma}
\label{lemma-hatX-decay} For almost every
$\omega\in\{0\in\scrC_{\infty,\alpha}\}$ and
every~$x\in\scrC_{\infty,\alpha}(\omega)$ there exists random variable
$C(\omega,x)<\infty$ such that
\begin{equation}
%\label{}
\hat{\cmss P}^n_\omega(x,y)\le\frac{C(\omega,x)}{n^{d/2}},\qquad n\ge1.
\end{equation}
\end{lemma}

We remark that the reversibility of the random walk, and the fact that $\pi_\omega\ge\alpha$ on~$\scrC_{\infty,\alpha}$, imply that  $\hat{\cmss P}^n_\omega(x,y)$ may also be bounded in terms of $C(\omega,y)$. Note that, unlike for $\cmss P$, the powers for which~$\hat{\cmss P}^n(x,y)$ is non-zero are not necessarily tied to the parity of~$y-x$.

\smallskip
Lemma~\ref{lemma-hatX-decay} will be implied by the fact that the Markov chain $\hat X$ obeys the ``usual'' $d$-di\-men\-sional isoperimetric inequality.  The connection between isoperimetric inequalities and heat-kernel decay
can be traced back to the work on elliptic PDEs done by Nash,
Moser and others. In its geometric form it was first proved using Sobolev
inequalities (Varopoulos~\cite{V85}). Alternative approaches use Nash inequalities
(Carlen, Kusoka and Stroock~\cite{CKS}), Faber-Krahn inequalities (Grigor'yan~\cite{G94}, Goel, Montenegro and Tetali~\cite{GM-Tetali}) and evolving sets (Morris and Peres~\cite{Morris-Peres}). 
The paper \cite{Morris-Peres} will serve us as a convenient reference. 

\smallskip
Consider a Markov chain on a countable state-space~$V$ with transition probability denoted by $\cmss P(x,y)$ and invariant measure denoted by~$\pi$. Define~$\cmss Q(x,y)=\pi(x)\cmss P(x,y)$ and for each~$S_1,S_2\subset V$, let
\begin{equation}
\label{QSS}
\cmss Q(S_1,S_2)=\sum_{x\in S_1}\sum_{y\in S_2}\cmss Q(x,y).
\end{equation}
For each~$S\subset V$ with~$\pi(S)\in(0,\infty)$ we define
\begin{equation}
\label{PhiS}
\Phi_S=\frac{\cmss Q(S,S^\cc)}{\pi(S)}
\end{equation}
and use it to define the isoperimetric profile
\begin{equation}
\label{Phi-inf}
\Phi(r)=\inf\bigl\{\Phi_S\colon \pi(S)\le r\bigr\}.
\end{equation}
(Here~$\pi(S)$ is the measure of~$S$.)
It is easy to check that we may restrict the infimum to sets~$S$ that are connected in the graph structure induced on~$V$ by $\cmss P$.

The following summarizes Theorem~2 of~\cite{Morris-Peres}: Suppose that~$\cmss P(x,x)\ge\gamma$ for some~$\gamma\in(0,\ffrac12]$ and all~$x\in V$. Let~$\epsilon>0$ and~$x,y\in V$. Then
\begin{equation}
\label{MP-bound}
\cmss P^n(x,y)\le\epsilon\pi(y)
\end{equation}
for all~$n$ such that
\begin{equation}
\label{LK-bound}
n\ge 1+\frac{(1-\gamma)^2}{\gamma^2}\int_{4[\pi(x)\wedge\pi(y)]}^{4/\epsilon}\frac4{u\Phi(u)^2}\,\textd u.
\end{equation}
Note that, to prove the ``usual'' decay $\cmss P^n(x,y)\le cn^{-d/2}$, it suffices to show that the bound $\Phi(r)\le c r^{-1/d}$ holds for~$r$ sufficiently large.

\smallskip
We will adapt this machinery to the following setting
\begin{equation}
V=\scrC_{\infty,\alpha}(\omega),\quad\cmss P=\hat{\cmss P}_\omega^2\quad\text{and}\quad\pi=\pi_\omega,
\end{equation}
with the objects in \twoeqref{QSS}{Phi-inf} denoted by~$\widehat{\cmss Q}_\omega$, $\widehat\Phi_S^{(\omega)}$ and~$\widehat\Phi_\omega(r)$. However, to estimate $\widehat\Phi_\omega(r)$ we will replace $\hat{\cmss P}_\omega$ by the Markov chain with transition probabilities
\begin{equation}
%\label{}
\tilde{\cmss P}_\omega(x,y)=\frac{\omega_{xy}\1_{\{\omega_{xy}\ge\alpha\}}}{\tilde\pi_\omega(x)},
\qquad|x-y|=1,
\end{equation}
i.e., the random walk on~$V=\scrC_{\infty,\alpha}$ that can only use edges physically present in the infinite cluster. The quantity
\begin{equation}
%\label{}
\tilde\pi_\omega(x)=\sum_{y}\omega_{xy}\1_{\{\omega_{xy}\ge\alpha\}}
\end{equation}
denotes the corresponding stationary measure. 
We will use $\widetilde{\cmss Q}_\omega$, $\widetilde\Phi_S^{(\omega)}$ and~$\widetilde\Phi_\omega(r)$ to denote the objects in \twoeqref{QSS}{Phi-inf} for this Markov chain. 

\begin{lemma}
\label{lemma-hat-tilde}
There exists a constant~$c>0$ depending only on~$d$ and~$\alpha$ such that for any finite set $\Lambda\subset\scrC_{\infty,\alpha}$,
\begin{equation}
\label{tildePhi-hatPhi}
\widehat\Phi_\Lambda^{(\omega)}\ge c \widetilde\Phi_\Lambda^{(\omega)}.
\end{equation}
\end{lemma}

\begin{proofsect}{Proof}
The stationary measures $\pi_\omega$ and $\tilde\pi_\omega$ compare via
\begin{equation}
\label{tildepi-pi}
\pi_\omega(x)\ge\tilde\pi_\omega(x)\ge\frac\alpha{2d}\pi_\omega(x)
\end{equation}
Restricting $\hat{\cmss P}_\omega^2(x,y)$ to transitions with~$T_2=2$ shows
\begin{equation}
%\label{}
\hat{\cmss P}_\omega^2(x,y)\ge\sum_z \frac{\omega_{xz}\1_{\{\omega_{xz}\ge\alpha\}}}{\pi_\omega(x)}\,
\frac{\omega_{zy}\1_{\{\omega_{zy}\ge\alpha\}}}{\pi_\omega(z)}
\ge(\ffrac \alpha{2d})^2\tilde{\cmss P}_\omega^2(x,y).
\end{equation}
It follows that \eqref{tildePhi-hatPhi} holds with $c=(\ffrac\alpha{2d})^3$.
\end{proofsect}

Our next step involves extraction of appropriate bounds on surface and volume terms. As the infimum in \eqref{Phi-inf} can always be restricted to connected subsets of the Markov graph, and since the Markov graph underlying the quantity $\widetilde\Phi_\omega(r)$ is just the infinite cluster~$\scrC_{\infty,\alpha}$, we can restrict our attention to subsets of~$\scrC_{\infty,\alpha}$ that are connected in the usual sense.

\begin{lemma}
\label{lemma-adapt}
Let~$\theta>0$, $d\ge2$ and let~$\alpha$ be as above. Then there exists a constant~$c>0$ and random variable $R_1=R_1(\omega)$ with~$\BbbP(R_1<\infty)=1$ such that for a.e.\ $\omega\in\{0\in\scrC_{\infty,\alpha}\}$ and all~$R\ge R_1(\omega)$ the following holds: For any connected~$\Lambda\subset\scrC_{\infty,\alpha}\cap[-R,R]^d$ with
\begin{equation}
%\label{}
\tilde\pi_\omega(\Lambda)\ge R^{\theta}
\end{equation}
we have
\begin{equation}
\label{Q-actual}
\widetilde{\cmss Q}_\omega(\Lambda,\scrC_{\infty,\alpha}\setminus\Lambda)\ge
c\pi_\omega(\Lambda)^{\frac{d-1}d}.
\end{equation}
\end{lemma}

\begin{proofsect}{Proof}
Since~$\scrC_{\infty,\alpha}$ has the law of the infinite bond-percolation cluster, we will infer this from isoperimetry for the percolation cluster; cf.\ Theorem~\ref{thm-isoperimetry}. Let~$\pomegaalpha\Lambda$ denote the set of edges in~$\scrC_{\infty,\alpha}$ with one endpoint in~$\Lambda$ and the other in $\scrC_{\infty,\alpha}\setminus\Lambda$. We claim that
\begin{equation}
\label{Q-bd}
\widetilde{\cmss Q}_\omega(\Lambda,\scrC_{\infty,\alpha}\setminus\Lambda)\ge \frac{\alpha^2}{2d}\,| \pomegaalpha\Lambda|
\end{equation}
and
\begin{equation}
\label{vol-bd}
\tilde\pi_\omega(\Lambda)\le 2d|\Lambda|.
\end{equation}
Since $\Lambda$ obeys the conditions of Theorem~\ref{thm-isoperimetry}, once~$R\gg1$, we have $| \pomegaalpha\Lambda|\ge c_2|\Lambda|^{\frac{d-1}d}$, cf.\ equation~\eqref{isoperimetry} in Sect.~\ref{sec5}. Then~\eqref{Q-actual} will follow from \twoeqref{Q-bd}{vol-bd}.

It remains to prove \twoeqref{Q-bd}{vol-bd}. The bound \eqref{vol-bd} is implied by~$\tilde\pi_\omega(x)\le2d$. For \eqref{Q-bd}, since~$\tilde{\cmss P}_\omega^2$ represents two steps of a random walk, we get a lower bound on~$\widetilde{\cmss Q}_\omega(\Lambda,\scrC_{\infty,\alpha}\setminus\Lambda)$ by picking a site~$x\in\Lambda$ which has a neighbor~$y\in\Lambda$ that has a neighbor~$z$ on the outer boundary of~$\Lambda$. The relevant contribution is bounded as
\begin{equation}
%\label{}
\tilde\pi_\omega(x)\tilde{\cmss P}_\omega^2(x,z)\ge\tilde\pi_\omega(x)\frac{\omega_{xy}}{\tilde\pi_\omega(x)}\frac{\omega_{yz}}{\tilde\pi_\omega(y)}\ge\frac{\alpha^2}{2d}.
\end{equation}
Once~$\Lambda$ has at least two elements, we can do this for~$(y,z)$ ranging over all bonds in~$\pomegaalpha\Lambda$, so summing over $(y,z)\in\pomegaalpha\Lambda$ we get~\eqref{Q-bd}.
\end{proofsect}

Now we are finally ready to estimate the decay of~$\hat{\cmss P}_\omega^n(x,y)$:

\begin{proofsect}{Proof of Lemma~\ref{lemma-hatX-decay}}
It clearly suffices to prove this for~$x=0$. Pick
$\theta\in(0,\ffrac12)$ and let~$R$ be the largest $\ell_\infty$-distance the walk~$X$ can go on~$\scrC_\infty$ by
time~$T_1+\dots+T_{2n}$, i.e., by the time~$\hat X$ makes~$2n$
steps. Lemma~\ref{prop-percolace} tells us that the
largest jump~$\hat X$ can make in a box of side length~$n^2$
is~$O(\log n)$, and so~$R=O(n\log n)$. As the walk will not leave the box~$[-R,R]^d$ by time~$n$, we may restrict the infimum defining~$\widehat\Phi_\omega(r)$ to sets $\Lambda$ entirely contained in~$[-R,R]^d$. (This can be formally achieved also by modifying the Markov chain ``outside'' $[-R,R]^d$.) Moreover, invoking \eqref{tildePhi-hatPhi} we can instead estimate $\widetilde\Phi_\omega(r)$ which allows us to restrict to~$\Lambda\subset\scrC_{\infty,\alpha}\cap[-R,R]^d$ that are connected in the usual graph structure on~$\scrC_{\infty,\alpha}$.

We will now derive a bound on~$\widetilde\Phi_\Lambda^{(\omega)}$ for connected~$\Lambda\subset\scrC_{\infty,\alpha}(\omega)\cap[-R,R]^d$. Henceforth~$c$ denotes a generic constant. 
If~$\pi_\omega(\Lambda)\ge R^{\theta}$, then \eqref{tildepi-pi} and \eqref{Q-actual} imply
\begin{equation}
%\label{}
\widetilde\Phi_\Lambda^{(\omega)}\ge c\pi_\omega(\Lambda)^{-\ffrac1d}.
\end{equation}
On the other hand, for~$\pi_\omega(\Lambda)<R^{\theta}$ the bound \eqref{Q-bd} yields
\begin{equation}
%\label{}
\widetilde\Phi_\Lambda^{(\omega)}\ge c\pi_\omega(\Lambda)^{-1}\ge cR^{-\theta}.
\end{equation}
From Lemma~\ref{lemma-hat-tilde} we conclude that
\begin{equation}
%\label{}
\widehat\Phi_\omega(r)\ge c\widetilde\Phi_\omega(r)\ge c (r^{-1/d}\wedge R^{-\theta})
\end{equation}
once~$R\ge R_1(\omega)$. The crossover
between the two regimes occurs when~$r=R^{d\theta}$ which (due to~$\theta<\ffrac12$) is much less than $4/\epsilon$ once $\epsilon\approx n^{-d/2}$. The relevant integral is thus bounded by
\begin{equation}
%\label{}
\int_{4[\pi(x)\wedge\pi(y)]}^{4/\epsilon}\frac4{u\widehat\Phi_\omega(u)^2}\,\textd u
\le c_1 R^{2\theta}\log R+c_2\epsilon^{-2/d}\le c_3\epsilon^{-2/d}
\end{equation}
for some constants~$c_1,c_2,c_3>0$. Setting~$\epsilon$ proportional to
$n^{-d/2}$ and noting $\gamma\ge(\ffrac\alpha{2d})^2$, the right-hand side
is less than~$n$ and $\cmss P^n(0,x)\le cn^{-d/2}$ for
each~$x\in\scrC_\infty\cap[-R,R]^d$. As $\cmss P^n(0,x)=0$ for~$x\not\in[-R,R]^d$, 
the bound holds in general. This proves the claim for even~$n$;
for odd~$n$ we just concatenate this with a single step of the
random walk.
\end{proofsect}

\subsection{Integral bound}
We now want to link the estimates on $\hat{\cmss P}$ to a heat-kernel type bound for the walk~$X$. Specifically, we will prove the following estimate:

\begin{proposition}
\label{prop-cond-bound} For almost every
$\omega\in\{0\in\scrC_{\infty,\alpha}\}$, there exists a constant
$C=C(\omega)<\infty$ such that for every $\ell\ge1$ and every
$n\ge1$,
\begin{equation}
\label{P-3.24}
P_{\omega,0}(\hat X_\ell=0,\,T_1+\dots+T_\ell\ge n)\le C(\omega)
\frac{\ell^{1-d/2}}n.
\end{equation}
and, in fact,
\begin{equation}
\label{eq:eqnum2}
\lim_{n\to\infty}\,n\, P_{\omega,0}(\hat
X_\ell=0,\,T_1+\dots+T_\ell\ge n)=0 \ \ \ \ \text{\rm a.s.}
\end{equation}
\end{proposition}

In order to prove this claim, we will need to occasionally refer to the Markov chain on environments ``from the point of view of the particle.'' Let~$\tau_x$ be the shift by~$x$ on $\Omega$ and let $\Omega_\alpha=\{0\in\scrC_{\infty,\alpha}\}$. We define a random shift~$\tau_{\hat X_1}\colon\Omega_\alpha\to\omega$ by sampling~$\hat X_1$ for the given~$\omega$ and applying~$\tau_x$ with~$x=\hat X_1$. This random map induces a Markov chain on $\Omega_\alpha^\Z$ via the iterated action of~$\tau_{\hat X_1}$. Define the measure
\begin{equation}
%\label{}
\Q_\alpha(\textd\omega)=Z\pi_\omega(0)\,
\BbbP(\textd\omega|0\in\scrC_{\infty,\alpha})
\end{equation}
where~$Z^{-1}=\E(\pi_\omega(0)|0\in\scrC_{\infty,\alpha})$. Let~$E_{\Q_\alpha}$ denote expectation with respect to~$\Q_\alpha$. We recall the following standard facts whose proof can be found in, e.g.,~\cite[Section~3]{BB}:

\begin{lemma}[Ergodicity of Markov chain on environments]
\label{lemma-ergodic}
The measure $\Q_\alpha$ is stationary and ergodic with respect to the Markov shift $\tau_{\hat X_1}$ on environments. In particular, if $f\in L^1(\Omega,\BbbP)$ then for $\Q_\alpha$-a.e.\ $\omega$ and for $P_{\omega,0}$-a.e.\ trajectory $\hat X=(\hat X_1,\hat X_2,\dots)$,
\begin{equation}
%\label{}
\lim_{\ell\to\infty}\frac1\ell\sum_{j=0}^{\ell-1}f(\tau_{\hat X_j}\omega)=E_{\Q_\alpha}(f).
\end{equation}
The convergence occurs also in~$L^1$ (i.e., under expectation~$E_{0,\omega}$ and, if desired, also~$E_{\Q_\alpha}$).
\end{lemma}

Recall our notation~$\FF_y$ for the finite component of~$\scrC_\infty\setminus\scrC_{\infty,\alpha}$ containing~$y$. For~$x\in\scrC_{\infty,\alpha}$, let
\begin{equation}
%\label{}
\GG_x'=\bigcup_{y\colon\omega_{xy}>0}\FF_y
\end{equation}
and let~$\GG_x$ denote the union of~$\GG_x'$ with all of its neighbors on~$\scrC_{\infty,\alpha}$. We will refer to this set as the \emph{weak component} incident to~$x$. Note that~$\GG_x$ is the set of vertices that can be visited by the walk~$X$ started at~$x$ by the time~$X$ steps again onto the strong component. 

\begin{lemma}
\label{lemma-exp-decay}
Recall that~$E_{\Q_\alpha}$ denotes expectation with respect to~$\Q_\alpha$ and let~$|\GG_x|$ be the number of sites in~$\GG_x$. Under the conditions of Lemma~\ref{prop-percolace}, we have $E_{\Q_\alpha}|\GG_0|<\infty$.
\end{lemma}

\begin{proofsect}{Proof}
This is an immediate consequence of \eqref{exponential}.
\end{proofsect}

Next we will estimate the expected time the random walk hides in such a component:

\begin{lemma}[Hiding time estimate]
\label{lemma-hiding}
Let~$d\ge2$ and set $c=4d\alpha^{-1}$. Then for all $x\in\Z^d$ and all $\omega$ such that~$x\in\scrC_{\infty,\alpha}$ and~$\GG_x$ is finite, we have
\begin{equation}
%\label{}
E_{\omega,x}(T_1)\le c|\GG_x|.
\end{equation}
\end{lemma}

\begin{proofsect}{Proof}
Fix $x\in\scrC_{\infty,\alpha}$ and let~$\GG_x$ be its incident weak component which we regard as a finite graph. Add a site~$\Delta$ to this graph and connect it by an edge to every site of~$\GG_x$ that has a strong bond to~$\scrC_{\infty,\alpha}\setminus\GG_x$. (Here~$\Delta$ represents the rest of~$\scrC_{\infty,\alpha}$; note that multiple edges between~$\Delta$ and sites of~$\GG_x$ are possible.) Equip each such edge with the corresponding conductance and call the resulting finite graph~$\HH_x$. Clearly, the random walk on~$\HH_x$ started at~$x$ and the corresponding random walk on~$\scrC_{\infty,\alpha}$ have the same law until they first hit~$\Delta$ (i.e., leave~$\GG_x$). In particular,~$T_1$ for the walk on~$\scrC_{\infty,\alpha}$ is stochastically dominated by $S_x$, the first time the walk on~$\HH_x$ returns back to its starting point.

Notice that $x\mapsto\pi_\omega(x)$ is an invariant measure of the walk on~$\HH_x$ provided we set
\begin{equation}
\pi_\omega(\Delta)=\sum_{x\in \GG_x}\sum_{y\in\scrC_{\infty,\alpha}\smallsetminus \GG_x}\omega_{xy}.
\end{equation}
Standard Markov chain theory tells us that~$z\mapsto(\tilde
E_zS_z)^{-1}$, where $\tilde E_z$ is the expectation with respect
to the walk on~$\HH_x$ started at~$z$, is an invariant
distribution and
\begin{equation}
\tilde E_x S_x=\frac{\pi_\omega(\HH_x)}{\pi_\omega(x)}.
\end{equation}
But~$x\in\scrC_{\infty,\alpha}$ implies that~$\pi_\omega(x)\ge\alpha$ while the bound~$\omega_{yz}\le1$ yields
\begin{equation}
\pi_\omega(\Delta)\le\pi_\omega(\GG_x)\le 2d|\GG_x|
\end{equation}
and
\begin{equation}
\pi_\omega(\HH_x)\le4d|\GG_x|.
\end{equation}
It follows that $E_{\omega,x}T_1\le\tilde E_xS_x\le(\ffrac{4d}\alpha)|\GG_x|$.
\end{proofsect}

\begin{proofsect}{Proof of Proposition~\ref{prop-cond-bound}}
For simplicity of the notation, let us assume that~$\ell$ is even; otherwise, replace all occurrences of~$\ffrac\ell2$ by~$\lceil\ffrac\ell2\rceil$. By reversibility of~$\hat X$, if~$k<\ell$,
\begin{equation}
%\label{}
P_{\omega,0}(\hat X_\ell=0,\,T_1+\dots+T_k\ge\ffrac n2)
=P_{\omega,0}(\hat X_\ell=0,\,T_\ell+\dots+T_{\ell-k+1}\ge\ffrac
n2).
\end{equation}
This means that the probability of interest is bounded by twice
the quantity on the left with~$k=\ffrac\ell2$. Chebyshev's inequality
then yields
\begin{equation}
\begin{aligned}
P_{\omega,0}(\hat X_\ell=0,\,T_1+\dots+T_\ell\ge n)
&\le 2P_{\omega,0}(\hat X_\ell=0,\,T_1+\dots+T_{\ell/2}\ge\ffrac n2)
\\
&\le\frac4nE_{\omega,0}\bigl(\1_{\{\hat X_\ell=0\}}(T_1+\dots+T_{\ell/2})\bigr).
\end{aligned}
\end{equation}
Conditioning on the position of~$\hat X$ at the times before and after~$T_j$ we then get
\begin{multline}
\label{3.25}
\qquad
P_{\omega,0}(\hat X_\ell=0,\,T_1+\dots+T_\ell\ge n)\\
\le\sum_{j=1}^{\ell/2}\sum_{x,y}\frac4n P_{\omega,0}(\hat X_{j-1}=x)\,E_{\omega,x}(T_1\,1_{\{\hat X_1=y\}})\,P_{\omega,y}(\hat X_{\ell-j}=0).
\qquad
\end{multline}
The calculation now proceeds by inserting uniform bounds for the last two terms on the right-hand side, and resumming the result using a stationarity argument.

Since~$\ell-j\ge\ell/2$, reversibility and Lemma~\ref{lemma-hatX-decay} tell us
\begin{equation}
\label{b1}
P_{\omega,y}(\hat X_{\ell-j}=0)=\frac{\pi_\omega(0)}{\pi_\omega(y)}
P_{\omega,0}(\hat X_{\ell-j}=y)\le\frac c{\ell^{d/2}}
\end{equation}
uniformly in $y\in\scrC_{\infty,\alpha}$ for some absolute constant~$c$. Furthermore, Lemma~\ref{lemma-hiding} gives
\begin{equation}
\sum_yE_{\omega,x}(T_1\,1_{\{\hat X_1=y\}})=E_{\omega,x}(T_1)\le c|\GG_x|
\end{equation}
where $\GG_x$ is the weak component incident to~$x$. Rewriting the sum over~$j$ as an ergodic average, Lemma~\ref{lemma-ergodic} with~$f(\omega)=|\GG_0|$ and Lemma~\ref{lemma-exp-decay} now show that, for all~$k\ge1$,
\begin{equation}
\label{b3}
\sum_{j=1}^k\sum_xP_{\omega,0}(\hat X_{j-1}=x)\,|\GG_x|=E_{\omega,0}\biggl(\,\sum_{j=0}^{k-1}|\GG_{\hat X_j}|\biggr)
\le C(\omega)k
\end{equation}
for a random constant~$C(\omega)$. Using \twoeqref{b1}{b3} in \eqref{3.25}, the desired bound \eqref{P-3.24} follows.

In order to prove the convergence to zero in \eqref{eq:eqnum2}, we note that
\begin{equation}
\sum_{n=1}^\infty P_{\omega,0}(\hat X_\ell=0,\,T_1+\dots+T_\ell\ge
n)=E_{\omega,0}\bigl(\1_{\{\hat X_\ell=0\}}\,(T_1+\dots+T_\ell)\bigr).
\end{equation}
The argument \twoeqref{3.25}{b3} shows that the expectation on the right is
finite a.s. Since $n\mapsto P_{\omega,0}(\hat X_\ell=0,\,T_1+\dots+T_\ell\ge n)$ is non-increasing, the claim follows by noting that, for any non-increasing non-negative sequence $(a_n)$ with $\limsup_{n\to\infty}na_n>0$, the sum $\sum_{n\ge1} a_n$ diverges.
\end{proofsect}

\subsection{Proof of the upper bound}
To turn \eqref{P-3.24} into the proof of Theorem~\ref{thm1}, we will also need the following standard fact from Markov chain theory:

\begin{lemma} \label{decreasing}
%\label{lemma}
The sequence $n\mapsto\cmss P_\omega^{2n}(0,0)$ is decreasing.
\end{lemma}

\begin{proofsect}{Proof}
Let $\langle f,g\rangle_\omega=\sum_{x\in\Z^d}\pi_\omega(x)f(x)g(x)$ denote a scalar product in $L^2(\Z^d,\pi_\omega)$. Then
\begin{equation}
%\label{}
\cmss P_\omega^{2n}(0,0) = \langle\delta_0,\cmss
P_\omega^{2n}\delta_0\rangle_\omega.
\end{equation}
Since $\cmss P_\omega$ is self-adjoint and $\Vert\cmss P_\omega\Vert_2\le1$, the sequence of operators $\cmss P_\omega^{2n}$ is decreasing.
\end{proofsect}

Now we put everything together and prove the desired heat-kernel upper bounds:

\begin{proofsect}{Proof of Theorem~\ref{thm1}(1)}
Introduce the random variable
\begin{equation}
%\label{}
R_n=\sup\{\ell\ge0\colon T_1+\dots+T_\ell\le n\}.
\end{equation}
The fact that $0\in\scrC_{\infty,\alpha}(\omega)$ yields
\begin{equation}
\sum_{m\ge n}P_{\omega,0}(X_m=0,\,R_m=\ell)
=P_{\omega,0}(\hat X_\ell=0,\,T_1+\dots+T_\ell\ge n).
\end{equation}
Proposition~\ref{prop-cond-bound} now implies
\begin{equation}
\label{eq:3.44}
\sum_{n\le m<2n}P_{\omega,0}(X_m=0,\,R_m=\ell)
\le C(\omega)\frac{\ell^{1-d/2}}n.
\end{equation}
By summing over $\ell=1,\dots,2n$ and using that $R_m\le 2n$ once $m\le 2n$ we derive
\begin{equation}
%\label{}
\sum_{n\le m<2n}\cmss P_\omega^m(0,0)
\le
\tilde C(\omega)\,\begin{cases}
n^{1-d/2},\qquad&d=2,3,
\\
n^{-1}\log n,\qquad&d=4,
\\
n^{-1},\qquad&d\ge5,
\end{cases}
\end{equation}
where $\tilde C$ is proportional to $C$. By Lemma \ref{decreasing}, $\cmss P_\omega^{2m}(0,0)$ 
is decreasing in~$m$ and so the sum on the left is bounded below by $\frac12n\cmss
P_\omega^{2n}(0,0)$. From here the claim follows.
\end{proofsect}

\begin{proofsect}{Proof of Theorem~\ref{thm1}(2)}
By \eqref{eq:eqnum2}, for each fixed~$\ell\ge1$ the sum in \eqref{eq:3.44} multiplied by~$n$ tends to zero as~$n\to\infty$. As~$\ell^{1-d/2}$ is summable in $d\ge5$, the uniform bound \eqref{eq:3.44} shows the same holds even under the sum over~$\ell\ge1$.
\end{proofsect}

\section{Examples with slow decay}
\label{sec4}\noindent
Here we provide proofs of Theorems~\ref{thm2} and~\ref{thm3}. The underlying ideas are very similar, but the proof of Theorem~\ref{thm2} is technically easier.

\subsection{Anomalous decay in $d\ge5$}
The proof of Theorem~\ref{thm2} will be based on the following
strategy: Suppose that in a box of side length~$\ell_n$ there
exists a configuration where a strong bond is separated from other
sites by bonds of strength~$\ffrac1n$, and (at least) one of these ``weak'' bonds is
connected to the origin by a ``strong'' path not leaving the box. Then
the probability that the walk is back to the origin at time~$n$ is
bounded below by the probability that the walk goes directly
towards the above pattern (this costs $\texte^{O(\ell_n)}$ of
probability) then crosses the weak bond (which costs~$\ffrac1n$),
spends time $n-2\ell_n$ on the strong bond (which costs only $O(1)$
of probability), then crosses a weak bond again (another factor
of~$\ffrac1n$) and then heads towards the origin to get there on
time (another~$\texte^{O(\ell_n)}$ term). The cost of this
strategy is $O(1)\texte^{O(\ell_n)}n^{-2}$ so if $\ell_n=o(\log
n)$ then we get leading order~$n^{-2}$.

\begin{proofsect}{Proof of Theorem~\ref{thm2}(1)}
Our task is to construct environments for which \eqref{lower-bd} holds.
For~$\kappa>\ffrac1d$ let $\epsilon>0$ be such that
$(1+4d\epsilon)/d<\kappa$. Let~$\B$ denote the set of edges in~$\Z^d$ and let~$\BbbP$ be an i.i.d.\ conductance
law on $\{2^{-N}\colon N\ge0\}^\B$ such that:
\begin{equation}
\label{p>pc}
\BbbP(\omega_b=1)>p_\cc(d)
\end{equation}
and
\begin{equation}
\BbbP(\omega_b=2^{-N})=cN^{-(1+\epsilon)},\qquad N\ge1,
\end{equation}
where $c=c(\epsilon)$ is adjusted so that the distribution is normalized. Let~$\hate_1$ denote the unit vector in the first coordinate direction. Define the scale
\begin{equation}
%\label{}
\ell_N=N^{(1+4d\epsilon)/d}
\end{equation}
and, given $x\in\Z^d$, let $A_N(x)$ be the event that the configuration near~$x$, $y=x+\hate_1$ and~$z=x+2\hate_1$ is as follows (see the comments before this proof):
%\settowidth{\leftmargini}{(11)}
\begin{enumerate}
\item[(1)]
$\omega_{yz}=1$ and $\omega_{xy}=2^{-N}$, while every other bond emanating out of~$y$ or~$z$ has $\omega_b\le 2^{-N}$.
\item[(2)]
$x$ is connected to the boundary of the box of side
length~$(\log\ell_N)^2$ centered at~$x$  by bonds with conductance
one.
\end{enumerate}
Since bonds with $\omega_b=1$ percolate and since $\BbbP(\omega_b\le 2^{-N})\sim N^{-\epsilon}$,  we have
\begin{equation}
%\label{}
\BbbP\bigl(A_N(x)\bigr)\ge cN^{-[1+(4d-2)\epsilon]}.
\end{equation}
Now consider a grid~$\G_N$ of sites
in~$[-\ell_N,\ell_N]^d\cap\Z^d$ that are spaced by
distance~$2(\log\ell_N)^2$. The events $\{A_N(x)\colon x\in\G_N\}$
are independent, so
\begin{equation}
\label{P-grid}
\BbbP\Bigl(\,\bigcap_{x\in\G_N}A_N(x)^\cc\Bigr)\le\exp\Bigl\{-c\Bigl(\frac{\ell_N}{(\log\ell_N)^2}\Bigr)^dN^{-[1+(4d-2)\epsilon]}\Bigr\}\le\texte^{-c N^{\epsilon}}
\end{equation}
and the intersection occurs only for finitely many~$N$.

 By the stretched-exponential decay of truncated
connectivities (Grimmett~\cite[Theorem~8.65]{Grimmett}), every connected
component of side length~$(\log\ell_N)^2$ in
$[-\ell_N,\ell_N]^d\cap\Z^d$ will eventually be connected to the
largest connected component in $[-2\ell_N,2\ell_N]^d\cap\Z^d$. We
conclude that there exists~$N_0=N_0(\omega)$
with~$\BbbP(N_0<\infty)=1$ such that once~$N\ge N_0$, the
event~$A_N(x)$ occurs for some even-parity site
$x=x_N(\omega)\in[-\ell_N,\ell_N]^d\cap\Z^d$ that is connected
to~$0$ by a path, $\text{Path}_N$, in $[-2\ell_N,2\ell_N]^d$, on which only the last $N_0$
edges---namely, those close to the origin---may have conductance
smaller than one.

We are now ready to employ the above strategy. Suppose~$N\ge N_0$
and let~$n$ be such that~$2^N\le 2n<2^{N+1}$. Let~$x_N$ be the
site in $[-\ell_N,\ell_N]^d\cap\Z^d$ for which~$A_N(x)$ occurs and
let~$r_N$ be the length of $\text{Path}_N$.
Let~$\alpha=\alpha(\omega)$ be the minimum
of~$\omega_b$ for~$b$ within~$N_0$ steps of the origin. The
passage from~$0$ to~$x_N$ in time~$r_N$ has probability at least
$\alpha^{N_0}(2d)^{-r_N}$, while staying on the bond~$(y,z)$ for
time~$2n-2r_N-2$ costs an amount which is bounded independently of~$\omega$. The
transitions across~$(x,y)$ cost order~$2^{-N}$ each. Hence we have
\begin{equation}
\label{P-4.6}
\cmss P_\omega^{2n}(0,0)\ge c\alpha^{2N_0}(2d)^{-2r_N}2^{-2N}.
\end{equation}
By the comparison of the graph-theoretic distance and the Euclidean distance (Antal and Pisztora~\cite{Antal-Pisztora}), we have $r_N\le c\ell_N$ once~$N$ is sufficiently large. Since~$n$ is of order~$2^N$ we are done.
\end{proofsect}

The argument for the second part follows very much the same strategy:

\begin{proofsect}{Proof of Theorem~\ref{thm2}(2)}
Let $(\lambda_n)$ be a sequence in the statement and suppose, without loss of generality, that $\lambda_1\gg1$. Let
\begin{equation}
q_n=\biggl(\frac12\frac{\log \lambda_n}{\log(2d)}\biggr)^{\ffrac14}
\end{equation}
and let $\{n_k\}$ be even numbers chosen as follows:
\begin{equation}
1-q_{n_1}^{-1}>p_\cc \quad\text{and}\quad
q_{n_{k+1}}>2q_{n_k}.
\end{equation}
Define an i.i.d.\ law $\BbbP$ on $(\{1\}\cup\{n_k\colon k\ge1\})^\B$ as follows:
\begin{equation}
\BbbP(\omega_b=1)=1-q_{n_1}^{-1}\quad\text{and}\quad
\BbbP(\omega_b=\ffrac1{n_k})=q_{n_k}^{-1}-q_{n_{k+1}}^{-1}.
\end{equation}
Let~$\eusm C_\infty$ denote the (a.s.\ unique) infinite connected component of edges with conductance one.

By following the argument in the proof of Theorem~\ref{thm2}(1), for almost every~$\omega$ and every $k$ large enough, we can find $x\in\eusm C_\infty$ such that:
%\settowidth{\leftmargini}{(11)}
\begin{enumerate}
\item
For $y=x+\hate_1$ and $z=x+2\hate_1$, we have $\omega_{y,z}=1$,
 and all other bonds emanating from~$y$ and~$z$ are of conductance~$1/n_k$.
\item
The chemical distance between $x$ and the closest point in~$\eusm C_\infty$ to the origin is less than~$q_{n_k}^{4}$.
\end{enumerate}
Explicitly, set $\ell_N=\theta q_{n_k}^4$ for some constant $\theta$ and let $A_n(x)$ be the event that~(1) holds and~$x$ is connected to the boundary of the box $x+[-(\log\ell_N)^2,(\log\ell_N)^2]^d$ by edges with strength one. Then $\BbbP(A_N(x))\ge cq_{n_k}^{-4d+2}=c\ell_N^{-d+\delta}$ for $\delta=(2d)^{-1}$. Plugging this in \eqref{P-grid} results in a sequence that is summable on~$k$ (note that $q_k$ increase exponentially). Percolation arguments, and the choice of~$\theta$, then ensure that (most of) the~$x$'s where~$A_N(x)$ occurs have a strong connection near the origin of length at most~$q_{n_k}^4$.

The argument leading to~\eqref{P-4.6}---with~$r_N$ replaced by~$q_{n_k}^{4}$---now gives
\begin{equation}
\cmss P_\omega^{n_k}(0,0)\geq
c\alpha^{2N_0}\frac{(2d)^{-2q_{n_k}^{4}}}{n_k^2}.
\end{equation}
By the choice of $q_n$, we are done.
\end{proofsect}

\subsection{Time-dependent environments}
Here we will prove Theorem~\ref{thm3}.
Let~$\BbbP$ be the Bernoulli measure on~$\B$ with parameter~$p>p_\cc(d)$. Let~$\eusm C_\infty$ denote the infinite component of occupied bonds. We define $\omega_b=1$ on occupied bonds and~$\omega_b=\ffrac1n$ on vacant bonds. The proof proceeds via three lemmas:

\begin{lemma}
\label{lemma-ac} Let~$Y=(Y_1,\dots,Y_n)$ be
the first~$n$ steps of the random walk on environment~$\omega$
conditioned to avoid bonds with~$\omega_b=\ffrac1n$. Let~$\tilde
X=(\tilde X_1,\dots,\tilde X_n)$ be the simple random walk on the percolation
cluster of~$\omega_b=1$. Then the the corresponding path measures are absolutely continuous with respect to each other and the Radon-Nikodym derivatives are (essentially) bounded away from zero and infinity, uniformly
in~$n$ and~$\omega\in\{0\in\eusm C_\infty\}$.
\end{lemma}

\begin{proofsect}{Proof}
Fix a sequence of sites $x_1,\dots,x_n\in\eusm C_\infty$ such
that~$\omega_{x_i,x_{i+1}}=1$ for all~$i=1,\dots,n-1$. Then the
probability that~$\tilde X$ executes this sequence is
$\prod_{i=0}^nd(x_i)^{-1}$, where~$d(x)$ is the degree of the
percolation cluster at~$x$. For~$Y$ we get
$C_n\prod_{i=0}^{n-1}\pi_\omega(x_i)^{-1}$, where~$C_n^{-1}$ is
the probability that the unconditioned random walk~$X$ has not
used a weak bond in its first~$n$-th steps. Since
\begin{equation}
%\label{}
\pi_\omega(x)-d(x)=O(\ffrac1n),
\end{equation}
the ratio of the products is bounded away from zero and infinity uniformly in~$n$ and the points $x_1,\dots,x_n$. But both path distributions are normalized and so~$C_n$ is bounded as well.
\end{proofsect}

Next we provide a lower bound on the probability that the walk~$X$ visits a given site in~$n$ steps. Let~$S_x$ be the first visit of~$X$ to~$x$,
\begin{equation}
%\label{}
S_x=\inf\{n\ge0\colon X_n=x\}.
\end{equation}
Then we have:

\begin{lemma}
\label{lemma-|x|}
For a.e.\ $\omega\in\{0\in\eusm C_\infty\}$ there is~$C=C(\omega)>0$ and a constant~$n_0<\infty$ such that for all $n\ge n_0$ and all~$x\in\eusm C_\infty$ satisfying $|x|\le\sqrt n$, we have
\begin{equation}
\label{tau-tail} P_{\omega,0}(S_x\le n)\ge
C(\omega)|x|^{-(d-2)}.
\end{equation}
\end{lemma}

\begin{proofsect}{Proof}
The choice of the conductance values ensures that the probability that~$X$ stays on~$\eusm C_\infty$ for
the first~$n$ steps is uniformly positive. Conditioning on this
event, and applying Lemma~\ref{lemma-ac}, it thus suffices to
prove \eqref{tau-tail} for the walk~$\tilde X$. The proof makes
use of Barlow's heat-kernel bounds for the random walk on
percolation cluster; cf~\cite[Theorem~1]{Barlow}.

Consider the continuous time version $\tilde X'$ of the
walk~$\tilde X$, i.e., $\tilde X'$ executes the same steps but at times that are i.i.d.\ exponential. 
By integrating the heat-kernel bounds we get that
the expected amount of time $\tilde X'$ spends at $x$ up to time
$n/2$ is at least $C(\omega)|x|^{-(d-2)}$. A similar calculation
shows that the expected time the walk $\tilde X'$  spends at $x$
\emph{conditioned on it hitting $x$} is uniformly bounded.
Therefore the probability of~$\tilde X'$ hitting $x$ before time
$n/2$ is at least $C(\omega)|x|^{-(d-2)}$. To get back to~$\tilde
X$, we need to subtract the probability that by continuous time
$n/2$ the walk $\tilde X'$ did more than $n$ discrete steps, which
is less than $\texte^{-cn}$. As~$|x|\le\sqrt n$, this cannot
compete with $|x|^{-(d-2)}$ once~$n$ is sufficiently large.
\end{proofsect}

We now define the notion of a \emph{trap} which is similar to that
underlying the event~$A_N(x)$ in the proof of Theorem~\ref{thm2}.
Explicitly, a trap is the triple of sites~$x,y,z$ with
$y=x+\hate_1$ and~$z=x+2\hate_1$ such that~$x\in\eusm C_\infty$
and such that all bonds emanating out of~$y$ and~$z$ are weak
except the bond between them. Let~$T(x)$ be the event that a trap
occurs at~$x$.

\begin{lemma}
\label{lemma-sum-|x|}
For a.e.\ $\omega\in\{0\in\eusm C_\infty\}$ there is $c<\infty$ and $n_1(\omega)<\infty$ such that
\begin{equation}
\sum_{\begin{subarray}{c}
x\colon|x|\le\sqrt n\\T(x)\text{\rm\ occurs}
\end{subarray}}
|x|^{-(2d-4)}\ge
\begin{cases}
c,\qquad&d\ge5,
\\
c\log n,\qquad&d=4,
\end{cases}
\end{equation}
for all $n\ge n_1$.
\end{lemma}

\begin{proofsect}{Proof}
This is a consequence of the Spatial Ergodic Theorem. Indeed,
let~$\Lambda_L=[-L,L]^d\cap\Z^d$ and note that the fraction
of~$\Lambda_L$ occupied by~$\{x\in\Lambda_L\colon T(x)\text{\rm\
occurs}\}$ converges a.s.\ to $\rho=\BbbP(T(0))>0$. But then also
the corresponding fraction in the annuli
$\Lambda_{2^{k+1}}\setminus\Lambda_{2^k}$ converges a.s.\
to~$\rho$. In particular, there is $k_0=k_0(\omega)$ such that
this fraction exceeds~$\rho/2$ for all $k\ge k_0$. Now take~$n$
and find~$k$ so that $2^k\le\sqrt n\le 2^{k+1}$. Bounding $|x|\le
2^{k+1}$ on the $k$-th annulus, we get
\begin{equation}
%\label{}
\sum_{\begin{subarray}{c}
x\colon|x|\le\sqrt n\\T(x)\text{\rm\ occurs}
\end{subarray}}
|x|^{-(2d-4)}
\ge\sum_{\ell=k_0}^k\frac\rho2\frac{|\Lambda_{2^{\ell+1}}\setminus\Lambda_{2^\ell}|}
    {(2^{\ell+1})^{2d-4}}.
\end{equation}
As $|\Lambda_{2^{\ell+1}}\setminus\Lambda_{2^\ell}|\ge (2^\ell)^d$, the result follows.
\end{proofsect}

We are now ready to prove the heat-kernel lower bounds~\eqref{eq2.7}:

\begin{proofsect}{Proof of Theorem~\ref{thm3}}
Pick~$\omega\in\{0\in\eusm C_\infty\}$ and let $x$ be a trap (i.e., event~$T(x)$ occurs and~$y$ and~$z$ are the endpoints of the ``trapped'' strong edge) with~$|x|<\frac14\sqrt n$. Let $U(x,k,\ell)$ be the event that the random walk starts at the origin, hits~$x$ for the first time at time~$k$, crosses the edge $(y,z)$, spends time $2n-k-\ell-2$ on this edge and then exits, and then arrives back to the origin in~$\ell$ units of time. Clearly,
\begin{equation}
%\label{}
P_{\omega,0}\bigl(U(x,k,\ell)\bigr)
\ge P_{\omega,0}(S_x=k)\,\frac cn\,\Bigl(1-\frac{\tilde c}n\Bigr)^{n-k-\ell-2}\,\frac cn\, P_{\omega,x}(S_0=\ell)
\end{equation}
where $c$ and $\tilde c$ are constants depending only on dimension. Reversibility tells us
\begin{equation}
%\label{}
P_{\omega,x}(S_0=\ell)=P_{\omega,0}(S_x=\ell)\frac{\pi_\omega(0)}{\pi_\omega(x)}\ge c
P_{\omega,0}(S_x=\ell)
\end{equation}
and so
\begin{equation}
%\label{}
P_{\omega,0}\bigl(U(x,k,\ell)\bigr)\ge c n^{-2}\,P_{\omega,0}(S_x=k)P_{\omega,0}(S_x=\ell).
\end{equation}
Denote
\begin{equation}
%\label{}
U(x)=\bigcup_{\begin{subarray}{c}
1\le k\le\ffrac n5\\\1\le \ell\le\ffrac n5
\end{subarray}}
U(x,k,\ell).
\end{equation}
Using the disjointness of~$U(x,k,\ell)$ for different~$k$ and~$\ell$ and invoking Lemma~\ref{lemma-|x|},
\begin{equation}
%\label{}
P_{\omega,0}\bigl(U(x)\bigr)\ge C(\omega)\,n^{-2}|x|^{-(2d-4)}.
\end{equation}
But, for~$n$ large enough, the events $\{U(x)\colon\,x\text{ is a trap}\}$ are disjoint because the restriction $k,\ell<\ffrac n5$ makes the walk spend more than half of its time at the strong bond constituting the trap. (This bond determines the trap entrance/exit point~$x$.) Hence,
\begin{equation}
%\label{}
\cmss P_\omega^{2n}(0,0)\ge
P_{\omega,0}\biggl(\,\bigcup_{x\colon|x|<\frac12\sqrt n}U(x)\biggr)\ge
C(\omega)\,n^{-2}\sum_{\begin{subarray}{c}
x\colon|x|\le\frac12\sqrt n\\T(x)\text{\rm\ occurs}
\end{subarray}}
|x|^{-(2d-4)}.
\end{equation}
Applying Lemma~\ref{lemma-sum-|x|}, the desired claim is proved.
\end{proofsect}

\section{Appendix: Isoperimetry on percolation cluster}
\label{sec5}\noindent
In this section we give a proof of isoperimetry of the percolation cluster which were needed in the proof of Lemma~\ref{lemma-adapt}. Consider bond percolation with parameter~$p$ and let~$\eusm C_\infty$ denote the a.s.\ unique infinite cluster. For~$\Lambda\subset\Z^d$ let~$\partial\Lambda$ denote the set of edges between~$\Lambda$ and~$\Z^d\setminus\Lambda$ and let $\pomega\Lambda$ denote those edges in~$\partial\Lambda$ that are occupied. We call~$\Lambda$ \emph{$\omega$-connected} if every two sites in~$\Lambda$ can be connected by a finite path that uses only the sites in~$\Lambda$ and whose every bond is occupied in~$\omega$. Then we have:

\begin{theorem}
\label{thm-isoperimetry}
For all $d\ge2$ and $p>p_\cc(d)$, there are positive and finite constants $c_1=c_1(d,p)$ and $c_2=c_2(d,p)$ and an a.s.\ finite random variable $R_0=R_0(\omega)$ such that for each $R\ge R_0$ and each $\omega$-connected~$\Lambda$ satisfying
\begin{equation}
%\label{}
\Lambda\subset\eusm C_\infty\cap[-R,R]^d\quad\text{and}\quad |\Lambda|\ge (c_1\log R)^{\frac d{d-1}}
\end{equation}
we have
\begin{equation}
\label{isoperimetry}
|\pomega\Lambda|\ge c_2|\Lambda|^{\frac{d-1}d}.
\end{equation}
\end{theorem}

This claim was the basic technical point of Benjamini and Mossel~\cite{Benjamini-Mossel} as well as of many subsequent studies of random walk on percolation cluster. Unfortunately, the proof of~\cite{Benjamini-Mossel} for the case~$d\ge3$ and~$p$ close to~$p_\cc(d)$ contains a gap. A different proof was recently given in Rau~\cite[Proposition~1.4]{Rau} but the argument is quite long and it builds (ideologically) upon a weaker version of \eqref{isoperimetry} proved by Mathieu and Remy~\cite{Mathieu-Remy}, whose proof is also rather long. Closely related estimates were derived in Barlow~\cite{Barlow}, but additional arguments are needed to extract~\eqref{isoperimetry}. 

For the convenience of the reader, and future reference, we provide a self-contained (and reasonably short) proof of Theorem~\ref{thm-isoperimetry} below. Our arguments are close to those of Benjamini and Mossel~\cite{Benjamini-Mossel} and they indicate that the seriousness of the gaps in~\cite{Benjamini-Mossel} has been somewhat exaggerated. An independent argument, based on exponential cluster repulsion, has simultaneously been found by Pete~\cite{Pete}.

\smallskip
Theorem~\ref{thm-isoperimetry} will be a consequence of the following, slightly more general estimate:

\begin{proposition}
\label{prop-bound}
For~$d\ge2$ and $p>p_\cc(d)$, there are~$c_2,c_3,\zeta\in(0,\infty)$ such that for all~$t>0$,
\begin{equation}
\label{bd}
\BbbP\bigl(\,\exists\Lambda\ni 0,\,\omega\text{\rm-connected},\,|\Lambda|\ge t^{\frac d{d-1}},\,|\pomega\Lambda|< c_2|\Lambda|^{\frac{d-1}d}\bigr)
\le c_3 \texte^{-\zeta t}.
\end{equation}
\end{proposition}

\begin{proofsect}{Proof of Theorem~\ref{thm-isoperimetry} from Proposition~\ref{prop-bound}}
Using translation invariance, the probability that there exists a set $\Lambda\subset\Z^d\cap[-R,R]^d$ with the properties listed in \eqref{bd} is bounded by a constant times~$R^d\texte^{-\zeta t}$. This applies, in particular, to sets~$\Lambda\subset\eusm C_\infty\cap[-R,R]^d$. Setting~$t=c_1\log R$ for~$c_1$ such that~$c_1\zeta>d+1$, this probability is summable on~$R$. By the Borel-Cantelli lemma, the corresponding event occurs only for finitely many~$R$. 
\end{proofsect}

The advantage of the formulation \eqref{bd} is that it links the tail bound on~$R_0$ to the cut-off on the size of~$|\Lambda|$. For instance, if we only care for~$|\Lambda|\ge R^\theta$ for some~$\theta\in(0,d)$, then~$\BbbP(R_0\ge R)$ decays exponentially with~$R^{\theta(1-1/d)}$.

\smallskip
As noted by Benjamini and Mossel~\cite{Benjamini-Mossel} the proof is quite straightforward in~$d=2$ and in any~$d$ once~$p$ is close to one.  However, to have a proof that works in $d\ge3$ all the way down to~$p_\cc$, we will have to invoke the ``static'' block-renormalization technique (Grimmett~\cite[Section~7.4]{Grimmett}). For each integer~$N\ge1$, consider the cubes
\begin{equation}
%\label{}
B_N(x)=x+\Z^d\cap[0,N]^d
\end{equation}
and
\begin{equation}
%\label{}
\tilde B_{3N}(x)=x+\Z^d\cap[-N,2N]^d
\end{equation}
Let~$G_N(x)$ be the event such that:
%\settowidth{\leftmargini}{(11)}
\begin{enumerate}
\item[(1)]
For each neighbor~$y$ of~$x$, the side of the block~$B_N(Ny)$ adjacent to~$B_N(Nx)$ is connected to the opposite side of~$B_N(Ny)$ by an occupied path.
\item[(2)]
Any two occupied paths connecting $B_N(Nx)$ to the boundary of $\tilde B_{3N}(Nx)$ are connected by an occupied path using only edges with both endpoints in~$\tilde B_{3N}(Nx)$.
\end{enumerate}
From Theorem~8.97 and Lemma~7.89 in Grimmett~\cite{Grimmett} we know that, for each~$p>p_\cc(d)$,
\begin{equation}
\label{good-to-one}
\BbbP\bigl(G_N(0)\bigr)\,\underset{N\to\infty}\longrightarrow\,1.
\end{equation}
By \cite[Theorem~7.65]{Grimmett}, for each~$p\in[0,1]$ there exists~$\eta_N(p)\in[0,1]$ with~$\eta_N(p)\uparrow1$ as $p\uparrow1$ such that the 0-1-valued process $\{\1_{G_N(x)}\colon x\in\Z^d\}$ is dominated from below by independent Bernoulli's with parameter~$\eta_N(p)$.

Given a finite set~$\Lambda\subset\Z^d$, let $\Lambda^{(N)}=\{x\in\Z^d\colon\Lambda\cap B_N(Nx)\ne\emptyset\}$
and define $\overline\Lambda^N$ to be the complement of the unique infinite component of~$\Z^d\setminus\Lambda^{(N)}$. We will also need a notation $\pstar\Delta$ for the inner site-boundary of a set $\Delta$,
\begin{equation}
\pstar\Delta=\bigl\{x\in\Delta\colon \exists y\in\Z^d\setminus\Delta\text{ with }|x-y|=1\bigr\},
\end{equation}
and $\diam\Lambda$ for the diameter of~$\Lambda$ in $\ell_\infty$-distance on $\Z^d$.
The crucial observation---which is where the setting of~\cite{Benjamini-Mossel} runs into a problem---is now as follows:

\begin{lemma}
\label{lemma-compare}
For~$\omega\in\Omega$, let $\Lambda\subset\Z^d$ be $\omega$-connected with $\overline\Lambda^N=\Delta$ and $\diam\Lambda\ge 3N$. If
\begin{equation}
\label{bd-small}
|\pomega\Lambda|<\frac1{2\cdot 3^d}|\pstar\Delta|
\end{equation}
then
\begin{equation}
\label{good-bd}
\bigl|\{x\in\pstar\Delta\colon G_N(x)^\cc\text{\rm\ occurs}\}\bigr|>\frac12|\pstar\Delta|.
\end{equation}
\end{lemma}

\begin{proofsect}{Proof}
Let $\Delta=\overline\Lambda^N$ and note that~$x\in\pstar\Delta$ implies $x\in\overline\Lambda^N$, i.e., $\Lambda\cap B_N(Nx)\ne\emptyset$. We claim that, for each~$x\in\pstar\Delta$,
\begin{equation}
G_N(x)\subset\bigl\{\tilde B_{3N}(Nx)\text{ contains an edge in }\pomega\Lambda\bigr\}.
\end{equation}
Indeed, if~$G_N(x)$ occurs then, by $\diam\Lambda\ge3N$, the box~$B_N(Nx)$ is connected to a site on the boundary of $\tilde B_{3N}(Nx)$ by an occupied path in~$\Lambda$. As~$x\in\pstar\Delta$ there exists a neighbor~$y\in\Delta^\cc$. Part (1) of the definition of~$G_N(x)$ ensures that there is another such path ``crossing'' $B_N(Ny)$; as $\Lambda\cap B_N(Ny)=\emptyset$, this path contains no sites in~$\Lambda$. By part~(2) of the definition of~$G_N(x)$, the two paths must be joined by an occupied path in~$\tilde B_{3N}(Nx)$ which then must contain an edge in~$\pomega\Lambda$.

Since each edge in~$\pomega\Lambda$ belongs to at most~$3^d$ distinct cubes~$\tilde B_{3N}(Nx)$ with~$x\in\pstar\Delta$, the number of boundary sites~$x\in\pstar\Delta$ where~$G_N(x)$ occurs is bounded by~$3^d|\pomega\Lambda|$, i.e.,
\begin{equation}
|\pstar\Delta|-\bigl|\{x\in\pstar\Delta\colon G_N(x)^\cc\text{\rm\ occurs}\}\bigr|\le 3^d|\pomega\Lambda|.
\end{equation}
Under the assumption \eqref{bd-small}, this implies \eqref{good-bd}.
\end{proofsect}

\begin{proofsect}{Proof of Proposition~\ref{prop-bound}}
Abbreviate~$c_4=(2\cdot 3^d)^{-1}$ and fix~$\Delta\subset\Z^d$ finite, connected with connected complement. Suppose~$\Lambda$ is $\omega$-connected with~$\overline\Lambda^N=\Delta$. Then
$|\Delta|\ge N^{-d}|\Lambda|$ and, invoking the standard isoperimetry on~$\Z^d$,
\begin{equation}
\label{iso-Delta}
|\pstar\Delta|\ge c_5|\Delta|^{\frac{d-1}d}\ge c_5N^{1-d}|\Lambda|^{\frac{d-1}d},
\end{equation}
where $c_5=c_5(d)>0$. Setting $c_2=c_4c_5N^{1-d}$ we then have
\begin{equation}
\label{include}
\bigl\{|\pomega\Lambda|< c_2|\Lambda|^{\frac{d-1}d}\bigr\}
\subset
\bigl\{|\pomega\Lambda|< c_4|\pstar\Delta|\bigr\}
\end{equation}
and also
\begin{equation}
\label{pDelta-t}
|\pstar\Delta|\ge c_5N^{1-d}t
\end{equation}
whenever $|\Lambda|\ge t^{\frac d{d-1}}$. We will suppose $t^{\frac d{d-1}}\ge(3N)^d$ to enable Lemma~\ref{lemma-compare}.

Equation \eqref{include}, Lemma~\ref{lemma-compare} and the fact that $\{\1_{G_N(x)}\colon x\in\Z^d\}$ stochastically dominates site percolation with parameter~$\eta_N(p)=1-\epsilon_N$ then yield
\begin{multline}
\label{perc-bound2}
\quad
\BbbP\bigl(\,\exists\Lambda\ni 0,\,\text{ $\omega$-connected},\,|\Lambda|\ge t^{\frac{d-1}d},\,\overline\Lambda^N=\Delta,\,|\pomega\Lambda|< c_2|\Lambda|^{\frac{d-1}d}\bigr)\\
\le\BbbP\left(\,\sum_{x\in\pstar\Delta}\1_{G_N(x)}\le\frac12|\pstar\Delta|\right)
\le2^{|\pstar\Delta|}(\epsilon_N)^{\frac12|\pstar\Delta|}.
\quad
\end{multline}
Here $2^{|\pstar\Delta|}$ bounds the number of possible subsets $\{x\in\pstar\Delta\colon G_N(x)^\cc\text{ occurs}\}$ of~$\pstar\Delta$.
To finish the proof, we need to sum over all eligible~$\Delta$'s. 

Let~$c_6=c_6(d)$ be a number such that~$c_6^n$ bounds the total number of connected sets~$\Delta\subset\Z^d$ with connected complement, containing the origin and having~$|\pstar\Delta|=n$. (The fact that this grows exponentially in~$n$ follows from the fact that~$\pstar\Delta$ is connected in an appropriate notion of adjacency on~$\Z^d$.) As~$\epsilon_N\to0$ by \eqref{good-to-one}, we can find~$N$ so that $2c_6\sqrt{\epsilon_N}\le\ffrac12$. Summing \eqref{perc-bound2} over all connected~$\Delta$ with connected complement that obey \eqref{pDelta-t} now gives
\begin{multline}
%\label{}
\qquad
\BbbP\bigl(\,\exists\Lambda\ni 0,\text{ $\omega$-connected},\,|\Lambda|\ge t^{\frac d{d-1}},\,|\pomega\Lambda|< c_2|\Lambda|^{\frac{d-1}d}\bigr)
\\
\le \sum_{n\ge c_5N^{1-d}t }2^{n}(\epsilon_N)^{\frac12n}c_6^n
\,\le\, \sum_{n\ge c_5N^{1-d}t }2^{-n}\le 2^{1-\lfloor c_5N^{1-d} t\rfloor},
\qquad
\end{multline}
where we also assumed that $2\epsilon_N^{1/2}\le1$ to get the first inequality.
Choosing the constants appropriately, this yields the desired claim.
\end{proofsect}

\section*{Acknowledgments}
\noindent The research of M.B.~was supported by the NSF
grant~DMS-0505356. The research of C.H.~was supported by the
NSF grant~DMS-0100445 and a grant from the University of
Washington Royalty Research Fund.
M.B.\ and G.K.\ would like to thank the Kavli
Institute for Theoretical Physics and their program ``Stochastic
Geometry and Field Theory: From Growth Phenomena to Disordered
Systems'' for support and hospitality that helped this project
take a definite shape. The KITP program was supported by the NSF under
the Grant No.~PHY99-07949. M.B.~also wishes to thank G\'abor Pete for discussions concerning the isoperimetry of the percolation cluster. Finally, we thank Martin Barlow and an anonymous referee for catching a number of inconsistencies in earlier versions of this manuscript.

\end{document}